\documentclass[11pt,reqno]{amsart}

\usepackage{amsmath,amsthm,amssymb,comment,fullpage}
\usepackage{braket}
\usepackage{mathtools}

\usepackage{caption}
\usepackage{times}
\usepackage{parskip}
\usepackage[T1]{fontenc}
\usepackage{mathrsfs}
\usepackage{latexsym}
\usepackage[dvips]{graphics}
\usepackage{epsfig}
\usepackage{amsmath,amsfonts,amsthm,amssymb,amscd}
\input amssym.def
\input amssym.tex
\usepackage{color}
\usepackage{hyperref}
\usepackage{url}
\newcommand{\bburl}[1]{\textcolor{blue}{\url{#1}}}

\usepackage{tikz}
\usepackage{tkz-tab}
\usepackage{tkz-graph}
\usetikzlibrary{decorations.pathreplacing}
\usetikzlibrary{shapes.geometric,positioning}

\newcommand{\burl}[1]{\textcolor{blue}{\url{#1}}}

\numberwithin{equation}{section}

\newtheorem{thm}{Theorem}[section]

\theoremstyle{plain}

\newtheorem{definition}[thm]{Definition}

\newtheorem{theorem}[thm]{Theorem}

\newcommand\be{\begin{equation}}
\newcommand\ee{\end{equation}}
\newcommand\bee{\begin{equation*}}
\newcommand\eee{\end{equation*}}
\newcommand\bea{\begin{eqnarray}}
\newcommand\eea{\end{eqnarray}}
\newcommand\beae{\begin{eqnarray*}}
\newcommand\eeae{\end{eqnarray*}}
\newcommand\bi{\begin{itemize}}
\newcommand\ei{\end{itemize}}
\newcommand\ben{\begin{enumerate}}
\newcommand\een{\end{enumerate}}
\newcommand\bc{\begin{center}}
\newcommand\ec{\end{center}}
\newcommand\ba{\begin{array}}
\newcommand\ea{\end{array}}

\newcommand\frakfamily{\usefont{U}{yfrak}{m}{n}}
\DeclareTextFontCommand{\textfrak}{\frakfamily}

\newcommand{\hr}[1]{\href{#1}{\url{#1}}}

\title{Counting on Euler and Bernoulli Number Identities}

\author{Arthur T. Benjamin}
\email{\textcolor{blue}{\href{mailto:benjamin@hmc.edu}{benjamin@hmc.edu}}}
\address{Department of Mathematics, Harvey Mudd College, Claremont, CA 91711}

\author{John Lentfer}
\email{\textcolor{blue}{\href{mailto:jlentfer@hmc.edu}{jlentfer@hmc.edu}}}
\address{Department of Mathematics, Harvey Mudd College, Claremont, CA 91711}

\author{Thomas C. Martinez}
\email{\textcolor{blue}{\href{mailto:tmartinez@hmc.edu}{tmartinez@hmc.edu}}}
\address{Department of Mathematics, Harvey Mudd College, Claremont, CA 91711}

\date{\today}

\begin{document}

\maketitle

\begin{abstract} While there are many identities involving the Euler and Bernoulli numbers, they are usually proved analytically or inductively. We prove two identities involving Euler and Bernoulli numbers with combinatorial reasoning via up-down permutations.
\end{abstract}

\section{Introduction}
Euler and Bernoulli numbers are numbers with interesting combinatorial properties, yet the combinatorial interpretation of the Bernoulli numbers has been little discussed so far, in part since they are not integers. In this paper, we discuss a combinatorial interpretation of identities that involve the Euler and Bernoulli numbers. As these numbers are sometimes indexed in different ways, we define then as follows:
\begin{definition}\label{def:def1}
Let $E_n$ be the $n$th Euler number and $B_n$ be the $n$th Bernoulli number; they are the unique sequence of numbers that satisfy 
\[
    \frac{2}{e^x+e^{-x}}\  =\  \sum_{n\,=\,0}^{\infty} \frac{E_n}{n!}\,x^n \quad \text{and} \quad \frac{x}{e^x-1}\ =\ \sum_{n\,=\,0}^{\infty} \frac{B_n}{n!}\,x^n.
\]
\end{definition}
In this paper, we show that the Euler numbers also satisfy
\begin{theorem}\label{thm:iden1} For $n\geq1$,
\begin{equation}\label{eqn:identity1}
    \sum_{j\,=\,0}^n\binom{2n}{2j}\,E_{2j}\ =\ 0.
\end{equation}
\end{theorem}
These numbers are also related to each other. We have the following identity that was previously proven analytically, as stated in \cite{Ent},
\begin{theorem}\label{thm:iden2} For $n\geq1$,
\begin{equation}\label{eqn:identity2}
    B_{2n}\ =\ \frac{2n}{2^{2n}(2^{2n}-1)}\,\sum_{j\,=\,0}^{n-1}\,\binom{2n-1}{2j}\,E_{2j}.
\end{equation}
\end{theorem}

Theorem \ref{thm:iden1}  has been proven combinatorially before, as seen in \cite{Men} and \cite{GS}, however our proof of Theorem \ref{thm:iden1} will set the stage for an original combinatorial proof of Theorem \ref{thm:iden2}. To prove these combinatorially, we need a combinatorial object to study.
Following \cite{Ent}, we will make use of up-down permutation, defined as follows:
\begin{definition}
For $n\geq 0$, an \emph{up-down permutation of length $n$} is a sequence of distinct numbers $a_1 a_2 \cdots a_n$ that satisfies  $a_1 < a_2 > a_3 < a_4 > \dots$. Let $U_n$ denote the number of up-down permutations that use all the numbers from $1$ through $n$. We allow the empty sequence to be counted so that $U_0 = 1$. 
\end{definition}
For example, $U_4 = 5$ counts the sequences 1324, 1423, 2314, 2413, and 3412. 
In \cite{Ent}, it was shown that $U_n$ is closely related to the Euler and Bernoulli numbers. Specifically, 
\begin{equation}\label{eqn:relation}
    U_{2n} = (-1)^n E_{2n} \quad \text{and} \quad  U_{2n-1} = \frac{2^{2n}(2^{2n}-1)}{2n}(-1)^{n-1}B_{2n}.
\end{equation}
Armed with these combinatorial objects, we prove Theorems \ref{thm:iden1} and \ref{thm:iden2}.

\section{Combinatorial Proofs}

\begin{proof}[Combinatorial Proof of Theorem \ref{thm:iden1}.] With equation \eqref{eqn:relation}, Theorem \ref{thm:iden1} can be rewritten as
$$ \sum_{j\,=\,0}^n\, \binom{2n}{2j} (-1)^j\, U_{2j}\ =\  0.$$
We prove this using the DIE method, as described in \cite{BQ}.

\textbf{Description:} In this step, we describe the unsigned objects. 
For $0\leq j\leq n$,  let $X$ denote an up-down permutation consisting of $2j$ elements from $\{1,2,\ldots,2n\}$. Then $X$ can be created in $\binom{2n}{2j}U_{2j}$ ways. Hence, the unsigned sum counts all up-down permutations of even length. 

\textbf{Involution:} Next we pair up the objects described above in such a way that paired objects have opposite sign in the sum. In other words, if object $X$ with length $2j$ is paired with $X^*$ of length $2j^*$, then $j$ and $j^*$ have opposite parity. 
Let $X\, =\, a_1\, a_2\, a_3\, a_4\, \cdots\, a_{2j}$ be an up-down permutation and let $S = \{s \in \{1,2,\dots, 2n\} \ |\ s \not\in X\}$ be the complementary subset. For example, if $2n = 8$, and $X = 6 \,8\, 2\, 7$, then $S = \{1, 3, 4, 5\}$.
For now, let's assume that $S$ and $X$ are nonempty. Let $y,z$ be the two largest elements of $S$, where $y<z$. We now consider two cases.

\underline{Case 1:} Suppose $z < a_1$. Then delete the first two elements of X. That is, let $X^*\, =\, a_3\, a_4\, \dots\, a_{2j}$ and $S^* = S  \cup \{a_1, a_2\}$. In our last example, $z = 5 < a_1 = 6$, so $X^* = 2 7$ and $S^* = \{1, 3, 4, 5, 6, 8\}$. 
Notice that $X^*$ is still an up-down permutation, and since $z < a_1 < a_2$, $a_1$ and $a_2$ are the largest elements of $S^*$. Notice that if we apply the involution to $X^*$, then Case 1 would not apply since $a_2 > a_3$ (in our example, $8 > 2$),  so the largest element of $S^*$ is greater than the first element of $X^*$. This leads us to the second case. 


\underline{Case 2:} Suppose $z > a_1$. Then append $y$ and $z$ to the beginning of $X$. That is, $X^*\, =\, y\, z\, a_1\, a_2\, a_3\, a_4\, \dots\, a_{2j}$, and $S^* = S - \{y, z\}$. For example, if $X = 2\, 7$, then $y = 6$ and $z = 8$, resulting in $X^* = 6\, 8\, 2\, 7$, as before. In general, since $y < z > a_1$, $X^*$ is guaranteed to be an up-down permutation, and if we apply the involution again, we will be in Case 1. 

Now we consider two edge cases that we previously ignored. If $S$ is empty, then we apply Case 1, so that $S^* = \{a_1, a_2\}$. If $X$ is empty, then $S = \{1, 2, \ldots, 2n\}$, and we apply Case 2 so that $S^* = \{1, 2, \ldots, 2n-2\}$ and $X^* =  (2n-1)\, 2n$. 

In summary, if $X$ is an up-down permutation of length $2j$, then our involution will either delete the first two elements of $X$, resulting in an up-down permutation of size $2j-2$ (and the deleted elements become the two largest missing elements) or it will 
append the two largest missing elements to the beginning of the $X$ resulting in an up-down permutation of size $2j+2$. Thus for every up-down permutation $X$, $(X^*)^* = X$. Moreover, $(-1)^{|X^*|/2} = (-1)^{j \pm 1}$, so our involution is sign-reversing.


%

\textbf{Exception:}
Since our involution is well-defined for every up-down permutation, including the empty permutation and those using all $2n$ elements, there are no exceptions to this rule. Every up-down permutation holds hands with an opposite-signed up-down permutation,
and the proof is complete. 
\end{proof}

\begin{proof}[Combinatorial Proof of Theorem \ref{thm:iden2}.]
With equation \eqref{eqn:relation}, Theorem \ref{thm:iden2} can be rewritten as
$$ \sum_{j\,=\,0}^{n-1} \,\binom{2n-1}{2j} \,(-1)^j\, U_{2j}\ =\ (-1)^{n-1}\,U_{2n-1}.$$
Once again, we prove this using the DIE method.

\textbf{Description:} As before, our unsigned summands count up-down permutations of even length. The only difference is that now the elements come from the odd-sized set  $\{1,2,\dots ,2n-1\}$. 

\textbf{Involution:} We use essentially the same involution as given in the previous proof. Let $X = a_1 a_2 \dots a_{2j}$, and let $y$ and $z$ denote the two largest unused elements, where $y < z$. If $z < a_1$, then we remove $a_1$ and $a_2$ from X. If $z > a_1$, we append $y z$ to the beginning of $X$. Note that $S$ will never be empty, as it will be of size $2n-1$ minus an even number, from the size of the permutation. If $X$ is empty, then $S = \{1, 2, \ldots, 2n-1\}$, and we apply Case 2 so that $S^* = \{1, 2, \ldots, 2n-3\}$ and $X^* =  (2n-2)\, (2n-1)$.

\textbf{Exception:} A potential problem arises when $j = n-1$. Here, $X = a_1 a_2 \dots a_{2n-2}$, so that $S=\{s\}$ consists of a single element. If $s < a_1$, then the involution works fine, since we can remove $a_1$ and $a_2$ from $X$ and these elements become the largest missing elements in the set $S^* = \{s, a_1, a_2\}$. But when $s > a_1$, then there is only one element of $S$, so we cannot apply two elements to $X$. For instance, if $n=3$ and $X\, =\, 1\, 5\, 3\, 4$, then $S = \{2\}$ contains a single element. Since $2 > 1$, our previous rule might suggest creating the permutation $2 1 5 3 4$, but this sequence is not an up-down permutation, nor does it have even length. On the other hand, it is a down-up permutation consisting of all the elements from the set $\{1, 2, 3, 4, 5\}$. In general, since $s > a_1$, the sequence $s a_1 a_2 \dots a_{2n-2}$ is a \textit{down-up} permutation, consisting of all the numbers from $\{1,2, \ldots, 2n-1\}$. But this allows us to count the exceptions in a nice way, since  $\{1, 2,\ldots, 2n-1\}$ has as many down-up permutations as up-down permutations. This can be seen through the bijection $f(a_1 a_2 \dots a_{2n-1}) = (2n - a_1) (2n - a_2) \dots (2n - a_{2n-1})$ which transforms a down-up permutation to an up-down permutation, and vice versa.

Thus the number of positions where the involution is undefined is $U_{2n-1}$. All of these occur when $j=n-1$, so they each have sign $(-1)^{n-1}$ in the summand, producing the right hand side of the identity. 
\end{proof}

\ 

\end{document}